\documentclass[a4paper]{article}
\usepackage{amsthm}
\usepackage[sumlimits]{amsmath}
\usepackage{amsfonts}
\DeclareMathOperator{\ad}{ad}
\newcommand{\bfi}{\bfseries\itshape}
\newtheorem{theorem}{Theorem}

\newtheorem{definition}[theorem]{Definition}
\newtheorem{remark}[theorem]{Remark}

\def\MM#1{\boldsymbol{#1}}

\newcommand{\pp}[2]{\frac{\partial #1}{\partial #2}}

\DeclareMathOperator{\Diff}{Diff}

\newcommand{\dede}[2]{\frac{\delta #1}{\delta #2}}
\newcommand{\dd}[2]{\frac{d #1}{d #2}}
\newcommand{\rd}{\mathrm{d}}

\begin{document}
\title{Multisymplectic formulation of fluid dynamics using the inverse
map} \author{C. J. Cotter, D. D. Holm and P. E. Hydon}
\label{firstpage}
\maketitle
\begin{abstract}{Multisymplectic systems, partial differential
    equations, fluid dynamics, conservation laws, potential vorticity}
  We construct multisymplectic formulations of fluid dynamics using
  the inverse of the Lagrangian path map. This inverse map -- the
  ``back-to-labels'' map -- gives the initial Lagrangian label of the
  fluid particle that currently occupies each Eulerian
  position. Explicitly enforcing the condition that the fluid
  particles carry their labels with the flow in Hamilton's principle
  leads to our multisymplectic formulation. We use the multisymplectic
  one-form to obtain conservation laws for energy, momentum and an
  infinite set of conservation laws arising from the
  particle-relabelling symmetry and leading to Kelvin's circulation
  theorem. We discuss how multisymplectic numerical
  integrators naturally arise in this approach.
\end{abstract}

\section{Introduction}

\begin{definition}
A system of partial differential equations (PDEs) is said to be
\emph{multisymplectic} if it is of the form
\[
K^{\alpha}_{ij}(\MM{z})z^j_{,\alpha} = \pp{H}{z^i},
\]
where each of the two-forms
\[
\kappa^\alpha = \frac{1}{2}K_{ij}^\alpha(\MM{z})\, \rd z^i\wedge\rd
z^j
\]
is closed. Here $\MM{z}$ is an ordered set of dependent variables,
total differentiation with respect to each independent variable
$q^\alpha$ is denoted by the subscript $\alpha$ after a comma, and the Einstein summation convention is used.
\end{definition}

The closed two-form $\kappa^\alpha$ is associated with the
independent variable $q^\alpha$; it is analogous to the symplectic
two-form for a Hamiltonian ordinary differential equation. Hence
there is a symplectic structure associated with each independent
variable. In the first of a series of papers, Bridges (1997)
pioneered the development of multisymplectic systems, showing that
the rich geometric structure that is endowed by the symplectic
two-forms can be used to understand the interaction and stability of
nonlinear waves. For many important PDEs, the multisymplectic
formulation has revealed hidden features that are important in
stability analysis. In order to preserve at least some of these
features in numerical simulations, Bridges \& Reich (2001)
introduced multisymplectic integrators, which generalise the
symplectic methods that have been widely used in numerical
Hamiltonian dynamics. Hydon (2005) showed that multisymplectic
systems of PDEs may be derived from Hamilton's principle whenever
the Lagrangian is affine in the first-order derivatives and contains
no higher-order derivatives. This can usually be achieved by
introducing auxiliary variables to eliminate the derivatives.

The aim of this paper is to provide a unified approach to producing
multisymplectic formulations of fluid dynamics, based on the inverse
map. Our approach covers all fluid dynamical equations that are
written in Euler-Poincar\'e form (Holm \emph{et al.}, 1998),
\emph{i.e.} all equations which arise due to the advection of fluid
material. First we use the inverse map to form a canonical
Euler-Lagrange equation (following the Clebsch representation given in
Holm \& Kupershmidt, 1983). Then the Lagrangian is made affine in the
space and time derivatives by using constraints that introduce
additional variables. Following Hydon (2005), we obtain a one-form
quasi-conservation law which, when it is pulled back to the space and
time coordinates, gives conservation laws for momentum and energy. We
also obtain a two-form conservation law that represents conservation
of symplecticity; when this is pulled back to the spatial coordinates,
it leads to a conservation law for vorticity.  The multisymplectic
version of Noether's Theorem yields an infinite space of conservation
laws from the particle-relabelling symmetry for fluid dynamics; these
conservation laws imply Kelvin's circulation theorem. The conserved
momentum that is canonically conjugate to the back-to-labels map plays
a key role in the derivation of the conservation laws. The
corresponding velocity is the convective velocity, whose geometric
properties are discussed in Holm \emph{et al.} (1986).

In this paper we show how the above constructions are made in
general, illustrating this with examples. We also discuss how
multisymplectic integrators can be constructed using these methods.
Sections \ref{review} and \ref{inverse map sec} review the relations
among multisymplectic structures, the Clebsch representation and the
momentum map associated with particle relabelling. Section
\ref{inverse map EPDiff} shows how to construct a multisymplectic
formulation of the Euler-Poincar\'e equation for the diffeomorphism
group (EPDiff), and derives the corresponding conservation laws,
including the infinite set of conservation laws that yield Kelvin's
circulation theorem. Section \ref{advected} extends this formulation
to the Euler-Poincar\'e equation with advected quantities. This is
illustrated by the incompressible Euler equation, showing how the
circulation theorem arises in the multisymplectic formulation.
Section \ref{numerics} sketches numerical issues in the
multisymplectic framework. Finally, Section \ref{summary} summarises
and outlines directions for future research.

\section{Review of multisymplectic structures}
\label{review}

This section reviews the formulation of multisymplectic systems and
their conservation laws, following Hydon (2005).

A system of partial differential equation (PDEs) is multisymplectic
provided that it can be represented as a variational problem with a
Lagrangian that is affine in the first derivatives of the dependent
variables:
\begin{equation}
\label{mslag}
L = L_j^\alpha(\MM{z})z^j_{,\alpha} - H(\MM{z}).
\end{equation}
The Euler-Lagrange equations are then
\begin{equation}
\label{Eul-Lag-eqns}K^{\alpha}_{ij}(\MM{z})z^j_{,\alpha} =
\pp{H}{z^i},
\end{equation}
where the functions
\begin{equation}
\label{Msymp-struct-matrix} K^{\alpha}_{ij}(\MM{z}) =
\pp{L^\alpha_j}{z^i}-\pp{L^\alpha_i}{z^j}
\end{equation}
are coefficients of the multisymplectic structure matrix. We define the
(closed) symplectic two-forms
\begin{equation}
\label{kappa} \kappa^\alpha = \frac{1}{2}K_{ij}^\alpha(\MM{z})\, \rd
z^i\wedge\rd z^j,
\end{equation}
and obtain the structural conservation law (Bridges, 1997).
\begin{equation}
\label{kappa law} \kappa^\alpha_{,\alpha} = 0.
\end{equation}
Hydon showed that the Poincar\'e Lemma leads to a one-form
quasi-conservation law
\begin{equation} (L^{\alpha}_jdz^j)_{,\alpha} =
\rd(L^\alpha_jz^j_{,\alpha}-H(\MM{z}))= \rd{L}, \label{ofcl}
\end{equation}
whose exterior derivative is (\ref{kappa law}).

Every one-parameter Lie group of point symmetries of the multisymplectic system (\ref{Eul-Lag-eqns}) is
generated by a differential operator of the form
\begin{equation}
\label{X}
X = Q^i(\MM{q},\MM{z})\pp{}{z^i} + (Q^i(\MM{q},\MM{z}))_{,\alpha}
\pp{}{z^i_{,\alpha}}.
\end{equation}
Noether's Theorem implies that if $X$ generates variational
symmetries, that is, if
\begin{equation}
\label{varsym}
XL = B^\alpha_{,\alpha}
\end{equation}
for some functions $B^\alpha$, then the
interior product of $X$ with the one-form quasi-conservation law
yields the conservation law
\begin{equation}
\label{noethm}
(L_j^{\alpha}Q^j-B^{\alpha})_{,\alpha}=0.
\end{equation}
This is the multisymplectic form of Noether's theorem.

Every multisymplectic system is invariant under translations in the
independent variables $\MM{q}$. For each of these symmetries,
Noether's theorem yields a conservation law
\[
(L_j^\alpha z^j_{,\beta}-L\delta^\alpha_\beta)_{,\alpha}=0.
\]
Such conservation laws can equally well be obtained by pulling back
the quasi-conservation law (\ref{ofcl}) to the base space of independent variables.
Commonly, the independent variables are spatial position $\MM{x}$ and time $t$.
Pulling back (\ref{ofcl}) to these base coordinates yields the
energy conservation law from the $\rd{t}$ component, and the
momentum conservation law from the remaining components. We shall
see the form of these conservation laws for fluid dynamics in later
sections.

\section{The inverse map and Clebsch representation}
\label{inverse map sec}

\subsection{Lagrangian fluid dynamics and the inverse map}

Lagrangian fluid dynamics provides evolution equations for particles
moving with a fluid flow. This is typically done by writing down a
flow map $\Phi$ from some reference configuration to the fluid
domain $\Omega$ at each instance in time. As the fluid particles
cannot cavitate, superimpose or jump, this map must be
a diffeomorphism.

For an $n$-dimensional fluid flow, the flow map
$\Phi:\,\mathbb{R}^n\times\mathbb{R}\mapsto\mathbb{R}^n$ given by
$\MM{x}=\Phi(\MM{l},t)$ specifies the spatial position  at time $t$
of the fluid particle that has \emph{label} $\MM{l}=\Phi(\MM{x},0)$.
The \emph{inverse map} $\Phi^{-1}$ gives the label of the particle
that occupies position $\MM{x}$ at time $t$ as the function
$\MM{l}=\Phi^{-1}(\MM{x},t)$.
The Eulerian velocity field $\MM{u}(\MM{x},t)$ gives the
velocity of the fluid particle that occupies position $\MM{x}$ at
time $t$ as follows:
\[
\MM{\dot{x}}(\MM{l},t)=\MM{u}(\MM{x}(\MM{l},t),t).
\]
Each label component $l_k(\MM{x},t)$
satisfies the advection law
\begin{equation}
\label{label eqn}
l_{k,t} + u_il_{k,i} = 0.
\end{equation}
Here ${}_{,t}$ and ${}_{,i}$ denote differentiation with respect to
$t$ and $x_i$ respectively. We use Cartesian coordinates and the
Euclidean inner product\footnote{This is only done for clarity and the
equations are easily extended to the case when the domain $\Omega$ is
a curved manifold.}, so we shall not generally distinguish between
`up' and `down' indices; summation from 1 to $n$ is implied whenever
an index is repeated.

\subsection{Clebsch representation using the inverse map}

A canonical variational principle for fluid dynamics may be formulated
by following the standard Clebsch procedure using the inverse map
(Seliger \& Whitham (1968), Holm \& Kupershmidt, 1983). The Clebsch
procedure begins with a functional $\ell[\MM{u}]$ of the Eulerian
fluid velocity $\MM{u}$, which is known as the \emph{reduced
  Lagrangian} in the context of Euler-Poincar\'e reduction (Holm
\emph{et al.}, 1998). One then enforces stationarity of the action
$S=\int\ell[\MM{u}]\rd t$ under the constraint that equation
(\ref{label eqn}) is satisfied by using a vector of $n$ Lagrange multipliers, which is denoted
as $\MM{\pi}$. These Lagrange multipliers are the conjugate
momenta to $\MM{l}$ in the course of the Legendre transformation to
the Hamiltonian formulation.  One may choose $\ell[\MM{u}]$ to be
solely the kinetic energy, which depends only on $\MM{u}$. More generally, $\ell$ will also depend on thermodynamic
Eulerian variables such as density, whose evolution may also be
accommodated by introducing constraints. These constraints are often
called the ``Lin constraints'' (Serrin, 1959). This idea was also
used in reformulating London's variational principle for superfluids
(Lin, 1963).
\begin{definition}[Clebsch variational principle using inverse map]
\label{clebsch}
The Clebsch variational principle using the inverse map is
\[
\delta \int_{t_0}^{t_1} \ell[\MM{u}] +
\int_{\Omega} \MM{\pi}\cdot(\MM{l}_t+\MM{u}\cdot\nabla\MM{l})\rd V(\MM{x})
\rd{t}=0,
\]
where $\MM{\pi}(\MM{x},t)$ are Lagrange multipliers which enforce the
constraint that particle labels $\MM{l}(\MM{x},t)$ are advected by the
flow.
\end{definition}
Taking the indicated variations leads to the following equations:
\begin{eqnarray}
\label{momentum map}
\delta \MM{u}:&&
\dede{\ell}{\MM{u}} + (\nabla\MM{l})^T\cdot\MM{\pi} = 0, \\
\nonumber
\delta \MM{\pi}:&&
\MM{l}_t+(\MM{u}\cdot\nabla)\MM{l} = 0, \\
\nonumber
\delta \MM{l}:&&
\MM{\pi}_t+\nabla\cdot(\MM{u}\MM{\pi}) = 0,
\end{eqnarray}
where
\begin{eqnarray}
\label{notation}
\left((\nabla\MM{l})^T\cdot\MM{\pi}\right)_i
:=
\pi_kl_{k,i},
\qquad
(\nabla\cdot(\MM{u}\MM{\pi}))_k
:=
(u_j\pi_k)_{,j},
\end{eqnarray}
and the variational derivative $\delta{\ell}/\delta{\MM{u}}$ is
defined by
\[
\ell[\MM{u}+\epsilon\MM{u}'] = \ell[\MM{u}] + \epsilon\int_{\Omega}
\dede{\ell}{\MM{u}}\cdot\MM{u}' \,\rd V(\MM{x}) +
\mathcal{O}(\epsilon^2) \,.
\]
\begin{remark}[Clebsch representation]
  In the language of fluid mechanics, the expression (\ref{momentum
    map}) for the spatial momentum $\MM{m}=\delta{\ell}/\delta\MM{u}$
  in terms of canonically conjugate variables $(\MM{l}, \MM{\pi})$ is
  an example of a ``Clebsch representation,'' which expresses the
  solution of the EPDiff equations (see below) in terms of canonical variables
  that evolve by standard canonical Hamilton equations.  This has been
  known in the case of fluid mechanics for more than 100 years.  For
  modern discussions of the Clebsch representation for ideal fluids,
  see, for example, Holm \& Kupershmidt (1983) and Marsden \&
  Weinstein (1983). In the language of geometric mechanics, the
  Clebsch representation is a momentum map.
\end{remark}
\subsection{Particle relabelling}

As the physics of fluids should be independent of the labelling of
particles, one may relabel the particles (by a diffeomorphism of the
flow domain) without changing the dynamics. This is called the
\emph{particle relabelling symmetry}; Noether's theorem applied to
this symmetry leads to the Kelvin circulation theorem. See Holm
\emph{et al.} (1998) for a modern description.

\subsection{Clebsch momentum map}
\begin{definition}
A {\bfi momentum map} is a map  $\mathbf{J}: T^{\ast}Q \rightarrow
\mathfrak{g}^\ast$ from the cotangent bundle $T^*Q$ of the
configuration manifold $Q$ to the dual $\mathfrak{g}^\ast$ of the
Lie algebra $\mathfrak{g}$ of a Lie group $G$ that acts on $Q$. The
momentum map is defined by the formula,
\begin{equation} \label{momentummapdef}
\mathbf{J} (\nu_q) \cdot \xi = \left\langle \nu_q, \xi_Q (q)
\right\rangle,
\end{equation}
where $\nu_q \in T ^{\ast} _q Q $ and $\xi \in \mathfrak{g}$. In
this formula $\xi _Q $ is the infinitesimal generator of the action
of ${G}$ on $Q$ associated with the Lie algebra element $\xi$, and
$\left\langle \nu_q, \xi_Q (q) \right\rangle$ is the natural pairing
of an element of $T ^{\ast}_q Q $ with an element of $T _q Q $.
\end{definition}

\begin{theorem}
  The Clebsch relation (\ref{momentum map}) defines a momentum map for
  the right action $\operatorname{Diff} (\Omega)$ of the
  diffeomorphisms of the domain $\Omega$ on the back-to-labels map
  $\MM{l}$.%
  \footnote{As we discuss later, this right action contrasts with
    fluid particle relabelling, which arises by the left action of the
    diffeomorphisms on the inverse map.}
\end{theorem}
\begin{proof}
  The spatial momentum in equation (\ref{momentum map}) may be
  rewritten as a map
  $\MM{J}_\Omega:\,T^*\Omega\mapsto\mathfrak{X}^*(\Omega)$ from the
  cotangent bundle of $\Omega$ to the dual $\mathfrak{X}^*(\Omega)$ of
  the vector fields $\mathfrak{X}(\Omega)$ given by
\begin{eqnarray}
\label{rightmommap}
\MM{J}_\Omega:\,\MM{m}\cdot \rd\MM{x}
 = - \Big( (\nabla\MM{l})^T\cdot\MM{\pi} \Big)\cdot \rd\MM{x}
 = - \,\MM{\pi} \cdot \rd\MM{l}
 =: - \,\pi_k \rd l_k
\,.
\end{eqnarray}
That is, $\MM{J}_\Omega$ maps the space of labels and their conjugate
momenta $(\MM{l},\MM{\pi})\in T^*\Omega$ to the space of one-form
densities $\MM{m}\in\mathfrak{X}^*(\Omega)$ on $\Omega$. The map
(\ref{rightmommap}) may be associated with the {\it right action}
$\MM{l}\cdot\eta$ of smooth invertible maps (diffeomorphisms) $\eta$
of the back-to-labels maps $\MM{l}$ by
composition of functions, as follows,
\begin{equation}\label{rightDiff}
\operatorname{Diff}(\Omega):\
\MM{l}\cdot\eta=\MM{l}\circ\eta
\,.
\end{equation}
The infinitesimal generator of this right action is obtained from its
definition, as
\begin{equation}\label{infgen-right}
X_{\Omega}(\MM{l})
:=
\frac{d}{ds}\Big|_{s=0}\Big(\MM{l}\circ\eta(s)\Big)
=
T\MM{l} \circ X
\,,
\end{equation}
in which the vector field $X \in \mathfrak{X}(\Omega)$ is tangent to
the curve of diffeomorphisms $\eta _s$ at the identity $s = 0$.
Thus, pairing the map $\MM{J}_\Omega$ with the vector field $X \in
\mathfrak{X}(\Omega)$ yields
\begin{align*}
\left\langle \MM{J}_\Omega (\MM{l}, \MM{\pi} ), X \right\rangle & =
-\,\langle\, \MM{\pi} \cdot \rd \MM{l} \,, X \,\rangle
\\& =
-\,\int_S
\pi_kl_{k,j}X_j (\MM{x})
\,\rd V(\MM{x})
\\& =
-\,\left\langle (\MM{l}, \MM{\pi}),
T\MM{l}\cdot X \right\rangle
\\& =
-\,\left\langle (\MM{l}, \MM{\pi}),
X_{\Omega}(\MM{l}) \right\rangle
\,,
\end{align*}
where $\left\langle \,\cdot\,, \,\cdot\,\right\rangle:\,T
^{\ast}_{\MM{l}}\Omega\times T_{\MM{l}}\Omega\mapsto\mathbb{R}$ is
the $L^2$ pairing of an element of $T ^{\ast}_{\MM{l}}\Omega $ (a
one-form density) with an element of $T_{\MM{l}}\Omega $ (a vector
field).
\\

\noindent
Consequently, the Clebsch map (\ref{momentum map}) satisfies the defining relation (\ref{momentummapdef}) to be a momentum map,
\begin{equation} \label{momentummap-JOmega}
\MM{J}(\MM{l}, \MM{\pi})
=
-\,\MM{\pi}
\cdot \rd \MM{l}
\,,
\end{equation}
with the $L^2$ pairing of the one-form density $-\,\MM{\pi}
\cdot \rd \MM{l}$ with the vector
field $X$.
\end{proof}

\begin{remark}
  Being the cotangent lift of the action of $\operatorname{Diff}
  (\Omega)$, the momentum map $\mathbf{J}_\Omega$ in
  (\ref{rightmommap}) is equivariant and Poisson. That is,
  substituting the canonical Poisson bracket into relation
  \eqref{momentummap-JOmega} yields the Lie-Poisson bracket on the
  space of $\MM{m}$'s. See, for example, Holm \& Kupershmidt (1983)
  and Marsden \& Weinstein (1983) for more explanation, discussion
  and applications. The momentum map property of the Clebsch
  representation guarantees that the canonically conjugate variables
  $(\MM{l},\MM{\pi})$ may be eliminated in favour of the spatial
  momentum $\MM{m}$. Before its momentum map property was understood,
  the use of the Clebsch representation to eliminate the canonical
  variables in favour of Eulerian fluid variables was a tantalising
  mystery (Seliger \& Whitham, 1968).

  Note that the right action of $\Diff(\Omega)$ on the inverse map is
  not a symmetry. In fact, as we shall see, the right action of
  $\Diff(\Omega)$ on the inverse map generates the fluid motion
  itself.
\end{remark}


\subsection{Elimination theorem}

Eliminating the canonically conjugate variables $(\MM{l}, \MM{\pi})$
produces an equation of motion for
$\MM{m}=\delta{\ell}/\delta\MM{u}$, which is constructed in the proof of the
following theorem:

\begin{theorem}[Elimination theorem]
The labels $\MM{l}$ and their conjugate momenta $\MM{\pi}$ may be
eliminated from the equations arising from the variational principle
(\ref{clebsch}) to obtain the weak form of the following equation of
motion for $\delta{\ell}/\delta\MM{u}$:
\[
\pp{}{t}\dede{\ell}{\MM{u}} +
\ad^*_{\MM{u}}\dede{\ell}{\MM{u}}
=0,
\]
where
\[
\ad^*_{\MM{u}}\MM{m} = \nabla\cdot(\MM{u}\MM{m}) + (\nabla\MM{u})^T\cdot
\MM{m}
\]
is defined by
\[
\langle \ad^*_{\MM{u}}\MM{m},\MM{w}\rangle = -\langle
\MM{m},\ad_{\MM{u}}\MM{w}\rangle
= \langle \MM{m},(\MM{u}\cdot\nabla)\MM{w} -
(\MM{w}\cdot\nabla)\MM{u}\rangle,
\]
and $\langle\cdot,\cdot\rangle$ is the $L^2$ inner-product.  This
equation of motion for ${\delta\ell}/{\delta\MM{u}}$ is the Euler-Poincar\'e
equation for the diffeomorphism group (EPDiff) (Holm \emph{et al.}, 1998).
\end{theorem}
\begin{proof}
Take the time-derivative of the inner product of ${\delta\ell}/{\delta\MM{u}}$
with a time-independent vector field $\MM{w}$:
\begin{eqnarray*}
\dd{}{t}\Bigg\langle \dede{\ell}{\MM{u}},\MM{w}\Bigg\rangle & = &
\dd{}{t}\Bigg\langle -(\nabla\MM{l})^T\cdot\MM{\pi}
, \MM{w} \Bigg\rangle =
\dd{}{t}\Bigg\langle -\MM{\pi},(\MM{w}\cdot\nabla)\MM{l}\Bigg\rangle, \\
& = & \Bigg\langle \nabla\cdot(\MM{u}\MM{\pi}),
(\MM{w}\cdot\nabla)\MM{l}\Bigg\rangle
+ \Bigg\langle \MM{\pi},
(\MM{w}\cdot\nabla)(\MM{u}\cdot\nabla)\MM{l}\Bigg\rangle
, \\
& = & \Bigg\langle \MM{\pi},
-(\MM{u}\cdot\nabla)(\MM{w}\cdot\nabla)\MM{l}+
(\MM{w}\cdot\nabla)(\MM{u}\cdot\nabla)\MM{l}\Bigg\rangle
= \Bigg\langle \MM{\pi},-\left(\ad_{\MM{u}}\MM{w}\cdot\nabla\right)\MM{l}
\Bigg\rangle , \\
&=& \Bigg\langle -(\nabla\MM{l})^T\cdot\MM{\pi},\ad_{\MM{u}}\MM{w}
\Bigg\rangle
 =  \Bigg\langle \dede{\ell}{\MM{u}},\ad_{\MM{u}}\MM{w}\Bigg\rangle
= -\,\Bigg\langle \ad^*_{\MM{u}}\dede{\ell}{\MM{u}},\MM{w}\Bigg\rangle,
\end{eqnarray*}
which is the (weak form of the) EPDiff equation.
\end{proof}

\subsection{Example: EPDiff($H^1$)}

To give a concrete example, consider EPDiff with $\ell[\MM{u}]$
being the $H^1_\lambda$-norm for $\MM{u}$. This is the
$n$-dimensional Camassa-Holm (CH) equation (Camassa \& Holm, 1993;
Holm \emph{et al.}, 1998; Holm \& Marsden, 2004), which has
applications in computational anatomy (Holm \emph{et al.}, 2004;
Miller \emph{et al.}, 2002). This system has the reduced Lagrangian
\[
\ell[\MM{u}] = \int_{\Omega}
\frac{1}{2}\,\left(|\MM{u}|^2+\lambda^2|\nabla\MM{u}|^2\right)\rd
V(\MM{x})  = \int_{\Omega}
\frac{1}{2}\,\left(u_iu_i+\lambda^2u_{i,j}u_{i,j}\right)\rd
V(\MM{x}) = \frac{1}{2}\|\MM{u}\|^2_{H^1_\lambda}.
\]
The EPDiff equation amounts to
\[
\pp{\MM{m}}{t} + (\MM{u}\cdot\nabla)\MM{m} +
(\nabla\MM{u})^T\cdot\MM{m} + \MM{m}\nabla\cdot\MM{u} = 0, \qquad
\MM{m} = (1-\lambda^2\nabla^2)\MM{u} \,.
\]
When $n=1$, these reduce to the Camassa-Holm (CH) equation,
\[
m_t + um_x + 2mu_x = 0,\qquad  m = u - \lambda^2u_{xx} \,.
\]

\subsection{Advected quantities}
To construct more general fluid equations we  shall include
advected quantities $a$ whose flow-rules are defined by
\begin{equation}
a_t + \mathcal{L}_{\MM{u}}a = 0,
\end{equation}
where $\mathcal{L}_{\MM{u}}$ is the Lie derivative. Such advected
variables typically arise in the potential energy or the thermodynamic
internal energy of an ideal fluid.  For example, advected scalars $s$
(as in salinity) satisfy
\begin{equation}
\pp{}{t}s + \mathcal{L}_{\MM{u}}s=0,
\quad \textrm{\emph{i.e.}} \quad
s_t + (\MM{u}\cdot\nabla)s = 0,
\end{equation}
and advected densities $\rho\,\rd V$ satisfy,
\begin{equation}
\pp{}{t}(\rho\,\rd V) + \mathcal{L}_{\MM{u}}(\rho\,\rd V)=0, \quad
\textrm{\emph{i.e.}} \quad \rho_t+\nabla\cdot(\rho\MM{u}) = 0. \,.
\label{advecD}
\end{equation}
A more extensive list of different types of advected quantity is
given in Holm \emph{et al.} (1998).

We write the reduced Lagrangian $\ell$ as a functional of the Eulerian
fluid variables $\MM{u}$ and $a$, and add further constraints to the
action $S$ to account for their advection relations,
\begin{equation}
\label{principle with advected qs}
S = \int \ell[\MM{u},a]\,\rd t
+ \int \rd t\int_{\Omega}\MM{\pi}\cdot(\MM{l}_t+
(\MM{u}\cdot\nabla)\MM{l}) + \phi(a_t+\mathcal{L}_{\MM{u}}a)
\,\rd V(\MM{x}).
\end{equation}
The Euler-Lagrange equations, which follow from the stationarity condition $\delta S=0$, are
\begin{eqnarray}
\delta\MM{u}:&&
\dede{\ell}{\MM{u}}  + (\nabla\MM{l})^T\cdot\MM{\pi}
 + \phi\diamond a = 0,
 \label{EL advected 1}\\
 \delta\MM{\pi}:&&
\MM{l}_t + (\MM{u}\cdot\nabla)\MM{l} = 0,
\nonumber\\
 \delta\MM{l}:&&
-\MM{\pi}_t - \nabla\cdot(\MM{u}\MM{\pi}) = 0,
\label{EL advected 3}\\
 \delta\phi:&&
a_t + \mathcal{L}_{\MM{u}}a = 0,
\nonumber \\
 \delta{a}:&&
-\phi_t-\mathcal{L}_{\MM{u}}\phi + \dede{\ell}{a} = 0,
\label{EL advected 5}
\end{eqnarray}
where the diamond operator ($\diamond$) is defined as the dual of the
Lie derivative operation $\mathcal{L}_{\MM{u}}$ with respect to the $L^2$ pairing. Explicitly, under integration by parts,
\begin{equation}
\int_\Omega (\phi\diamond a)\cdot\MM{u}\,\rd{V}(\MM{x})
= -\int_\Omega (\phi\mathcal{L}_{\MM{u}}a)\,\rd{V}(\MM{x}).
\end{equation}
\begin{remark}The map to the spatial momentum in equation (\ref{EL advected 1})
\begin{equation}
\dede{\ell}{\MM{u}} =: \MM{m}
= -\,\pi_A\nabla l^A -\, \phi\diamond a
\,,
\end{equation}
is again a momentum map, this time for the semidirect-product action
of the diffeomorphisms on $\Omega\times V^*$. Again the momentum map
property allows the canonical variables to be eliminated in favour of
the Eulerian quantities. As a result, eliminating the variables $\MM{l}$,
$\MM{\pi}$ and $\phi$ leads to the Euler-Poincar\'e equation with
advected quantities $a$.
\end{remark}

\begin{theorem}[Elimination theorem with advected quantities]
  The labels $\MM{l}$, their conjugate momenta $\MM{\pi}$ and the
  conjugate momentum $(\phi)$ to the advected quantities $(a)$ may be
  eliminated from equations (\ref{EL advected 1}-\ref{EL advected 5})
  to obtain the weak form of the Euler-Poincar\'e equation with
  advected quantities:
\[
\pp{}{t}\dede{\ell}{\MM{u}} + \ad^*_{\MM{u}}\dede{\ell}{\MM{u}}
= a\diamond\dede{\ell}{a}, \qquad a_t + \mathcal{L}_{\MM{u}}a = 0.
\]
\end{theorem}
\begin{proof}
  Take the time-derivative of the inner product of
  ${\delta\ell}/{\delta\MM{u}}$ with a function of $\MM{w}$:
\begin{eqnarray*}
  \dd{}{t}\Bigg\langle \dede{\ell}{\MM{u}},\MM{w}\Bigg\rangle & = &
  \dd{}{t}\Bigg\langle -(\nabla\MM{l})^T\cdot\MM{\pi} - \phi\diamond a
  , \MM{w} \Bigg\rangle
  =  \dd{}{t}\Bigg\langle -\MM{\pi},(\MM{w}\cdot\nabla)\MM{l}\Bigg\rangle
  + \dd{}{t}\Bigg\langle \phi,\mathcal{L}_{\MM{w}}a\Bigg\rangle \\
  & = & \Bigg\langle \nabla\cdot(\MM{u}\MM{\pi}),
  (\MM{w}\cdot\nabla)\MM{l}\Bigg\rangle
  + \Bigg\langle \MM{\pi},
  (\MM{w}\cdot\nabla)(\MM{u}\cdot\nabla)\MM{l}\Bigg\rangle
  + \Bigg\langle \dede{\ell}{a}\diamond a,\MM{w}\Bigg\rangle\\
  & & \quad + \Bigg\langle -\dede{\ell}{a}-\mathcal{L}_{\MM{u}}\phi,
  \mathcal{L}_{\MM{w}}a\Bigg\rangle
  + \Bigg\langle \phi,-\mathcal{L}_{\MM{w}}\mathcal{L}_{\MM{u}}a\Bigg\rangle
  , \\
  & = & \Bigg\langle \MM{\pi},
  -(\MM{u}\cdot\nabla)(\MM{w}\cdot\nabla)\MM{l}+
  (\MM{w}\cdot\nabla)(\MM{u}\cdot\nabla)\MM{l}\Bigg\rangle
  + \Bigg\langle \dede{\ell}{a}\diamond a,\MM{w}\Bigg\rangle
  \\
  & & \quad
  +\Bigg\langle\phi,\mathcal{L}_{\MM{u}}\mathcal{L}_{\MM{w}}a
  -\mathcal{L}_{\MM{w}}\mathcal{L}_{\MM{u}}a \Bigg\rangle
  , \\
  & =&\Bigg\langle \MM{\pi},-\left(\ad_{\MM{u}}\MM{w}\cdot\nabla\right)\MM{l}
  \Bigg\rangle
  + \Bigg\langle \phi,\mathcal{L}_{\ad_{\MM{u}}\MM{w}}a
  \Bigg\rangle + \Bigg\langle \dede{\ell}{a}\diamond a,\MM{w}\Bigg\rangle, \\
  &=& \Bigg\langle -(\nabla\MM{l})^T\cdot\MM{\pi}
  -\phi\diamond a,\ad_{\MM{u}}\MM{w}
  \Bigg\rangle
  + \Bigg\langle \dede{\ell}{a}\diamond a,\MM{w}\Bigg\rangle , \\
  & = & \Bigg\langle \dede{\ell}{\MM{u}},\ad_{\MM{u}}\MM{w}\Bigg\rangle
  + \Bigg\langle \dede{\ell}{a}\diamond a,\MM{w}\Bigg\rangle
  =  \Bigg\langle -\ad^*_{\MM{u}}\dede{\ell}{\MM{u}}
  + \dede{\ell}{a}\diamond a,\MM{w}\Bigg\rangle.
\end{eqnarray*}
\end{proof}

\begin{remark}These Euler-Poincar\'e equations with advected
  quantities cover all conservative fluid equations which describe the
  advection of material. For a large collection of examples, see Holm
  \emph{et al.} (1998).
\end{remark}

\subsection{Example: Incompressible Euler equations}

As an example, consider the reduced Lagrangian for the incompressible
Euler equations
\[
\ell[\MM{u},\rho,p] = \int_\Omega \frac{\rho|\MM{u}|^2}{2} +
p(1-\rho)\,\rd V(\MM{x}) \,.
\]
Here $\rho(\MM{x},t)$ is the ratio of the local fluid density to the
average density over $\Omega$; this is governed by the continuity
equation (\ref{advecD}). The pressure $p$ is a Lagrange multiplier
that fixes the incompressibility constraint $\rho=1$. The
variational derivatives in this case are
\[
\dede{\ell}{\MM{u}} = \rho\MM{u}, \qquad \dede{\ell}{\rho} =
\frac{|\MM{u}|^2}{2}-p,\qquad \dede{\ell}{p} =1-\rho\,.
\]
Consequently, the Euler-Poincar\'e equations become
\begin{eqnarray*}
(\rho\MM{u})_t + (\MM{u}\cdot\nabla)(\rho\MM{u}) +
\rho\MM{u}(\nabla\cdot\MM{u}) + \rho(\nabla\MM{u})^T\cdot \MM{u} &=&
\rho\nabla\left(\frac{|\MM{u}|^2}{2} -p\right), \\
\rho_t+\nabla\cdot(\rho\MM{u}) &=& 0, \\ \qquad \rho&=&1,
\end{eqnarray*}
and rearrangement gives the Euler fluid equations,
\[
\MM{u}_t + (\MM{u}\cdot\nabla)\MM{u} = -\nabla p, \qquad \nabla\cdot\MM{u}=0.
\]
\section{Inverse map multisymplectic formulation for EPDiff($H^1$)}
\label{inverse map EPDiff} As we now have a canonical variational
principle for fluid dynamics \emph{via} the inverse map, one may
obtain its multisymplectic formulation by extending the phase space
so that the Lagrangian is affine in the space and time derivatives.
In this section we show how to do this for EPDiff($H^1$) as
discussed in the previous section.

\subsection{Affine Lagrangian for EPDiff($H^1$)}

After introducing the inverse map constraint, the Lagrangian
becomes
\[
L = \frac{1}{2}u_iu_i +
\frac{\lambda^2}{2}u_{i,j}u_{i,j} +
\pi_k\left(l_{k,t} + u_jl_{k,j}\right).
\]
Any high-order
derivatives and nonlinear functions of first-order derivatives must
now be removed from the Lagrangian to make it affine. We introduce a
tensor variable
\[
W_{ij} = u_{i,j}\,;
\]
this relationship may be enforced by using Lagrange multipliers.
However, it turns out that the multipliers can be eliminated and the
Lagrangian becomes
\begin{equation}
\label{epmslag}
L = \frac{1}{2}u_iu_i - \frac{\lambda^2}{2}
W_{ij}W_{ij} +\lambda^2W_{ij}u_{i,j}+ \pi_k\left(l_{k,t} +
u_jl_{k,j}\right),
\end{equation}
which is now affine in the space and time derivatives of $\MM{u}$,
$W$, $\MM{l}$ and $\MM{\pi}$.

\subsection{Multisymplectic structure}
The Euler-Lagrange equations for the affine Lagrangian (\ref{epmslag}) are
\begin{eqnarray*}
\delta u_i:&& u_i  - \lambda^2W_{ij,j}+ \pi_k l_{k,i} = 0,
\\
 \delta l_k:&&
-\pi_{k,t} - (\pi_k u_j)_{,j} = 0.
\\
 \delta \pi_k:&&
l_{k,t} + u_j l_{k,j} = 0,
\\
 \delta W_{ij}:&&
-\lambda^2 W_{ij} + \lambda^2 u_{i,j}
 = 0 .
\end{eqnarray*}
These equations possess the following multisymplectic structure as
in equation (\ref{Eul-Lag-eqns}):
\[
\begin{pmatrix}
0 & \pi_k\partial_i &  & -\lambda^2\partial_j \\
-\pi_k\partial_i & 0 & -\partial_t-u_j\partial_j & 0 \\
0 & \partial_t+u_j\partial_j & 0 & 0 \\
\lambda^2\partial_j & 0 & 0 & 0 \\
\end{pmatrix}
\begin{pmatrix}
u_i \\
l_k \\
\pi_k \\
W_{ij}\\
\end{pmatrix}
= \nabla H,
\]
where $\partial_t=\partial/\partial t,\ \partial_i=\partial/\partial x_i$, and
\[
H = -\left(\frac{1}{2}u_iu_i -\frac{\lambda^2}{2}W_{ij}W_{ij}\right) = -\left(\frac{1}{2}|\MM{u}|^2 -\frac{\lambda^2}{2}|W|^2\right).
\]

\subsection{One-form quasi-conservation law}

For our multisymplectic formulation of EPDiff($H^1$), the independent
variables are
\[
q^j = x_j, \quad j=1,\ldots n, \qquad q^{n+1} = t,
\]
and the dependent variables are
\begin{eqnarray}
z^i = u_i, \qquad z^{n + k} =
l_k, \qquad z^{2n + k} = \pi_k,\qquad z^{(i+2)n+j}=W_{ij}\,,
\end{eqnarray}
where $i,j$ and $k$ range from $1$ to $n$. Comparing (\ref{epmslag})
with (\ref{mslag}) gives the following non-zero components
$L^{\alpha}_j$:
\[
L^j_i = \lambda^2W_{ij}, \qquad L_{n+k}^j = \pi_ku_j, \qquad
L_{n+k}^{n+1} = \pi_k, \qquad i,j,k=1,\ldots,n.
\]
Therefore the one-form quasi-conservation law amounts to
\begin{equation}
\label{epqcl}
\left(\pi_k\rd l_k\right)_{,t}+\left(\lambda^2W_{ij}\rd u_i+\pi_ku_j\rd l_k\right)_{,j}=\rd L.
\end{equation}
The exterior derivative of this expression yields the structural
conservation law
\begin{equation}
\label{epsympcl}
\left(\rd\pi_k\wedge\rd l_k\right)_{,t}+\left(\lambda^2\rd W_{ij}\wedge\rd u_i+u_j\rd\pi_k\wedge\rd l_k+\pi_k\rd u_j\wedge\rd l_k\right)_{,j}=0.
\end{equation}

\subsection{Conservation of energy}
For EPDiff($H^1$), the
$\rd t$-component of the pullback of the one-form conservation law (\ref{epqcl}) gives
\[
\left(\pi_kl_{k,t} - L\right)_{,t} +\left(
\lambda^2W_{ij}u_{i,t} + \pi_ku_jl_{k,t} \right)_{,j} = 0.
\]
In terms of $\MM{u}$ and its derivatives, this amounts to
\[
\left(u_im_i-\frac{1}{2}u_iu_i-\frac{\lambda^2}{2}u_{i,j}u_{i,j}\right)_{,t}
+\left(\lambda^2u_{i,j}u_{i,t}+u_iu_jm_i\right)_{,j} =
0,
\]
where
\[
m_i=u_i-\lambda^2u_{i,kk}.
\]
This is the energy conservation law for EPDiff($H^1$).

\subsection{Conservation of momentum}
Similarly, the conservation law that is associated with translations in the $x_i$-direction is
\[
\left(\pi_kl_{k,i}\right)_{,t} +\left( \lambda^2
W_{kj}u_{k,i} + \pi_ku_jl_{k,i} - \delta_{ij}L
\right)_{,j}=0,
\]
which amounts to the momentum conservation law
\[
m_{i,t}+\left(\lambda^2u_{k,i}u_{k,j}-u_jm_i-\delta_{ij}\left(\frac{1}{2}u_ku_k+\frac{\lambda^2}{2}u_{k,l}u_{k,l}\right)\right)_{,j}.
\]

\subsection{Conservation of vorticity}
Next, consider the coefficient of each $\rd x_r\wedge\rd x_s$ in the
pull-back of the structural (two-form) conservation law
(\ref{epsympcl}). This is
\[
\left(\pi_{k,r}l_{k,s}-\pi_{k,s}l_{k,r}\right)_{,t}
+\left(\lambda^2\left(W_{ij,r}u_{i,s}-W_{ij,s}u_{i,r}\right)
+u_j\left(\pi_{k,r}l_{k,s}-\pi_{k,s}l_{k,r}\right)
+\pi_k\left(u_{j,r}l_{k,s}-u_{j,s}l_{k,r}\right)\right)_{,j}=0,
\]
which amounts to
\[
\left(m_{r,s}-m_{s,r}\right)_{,t}+\left(\lambda^2\left(u_{i,s}u_{i,jr}-u_{i,r}u_{i,js}\right)+(u_jm_r)_{,s}-(u_jm_s)_{,r}\right)_{,j}=0.
\]
One can regard this as a vorticity conservation law for
EPDiff($H^1$); it is a differential consequence of the momentum
conservation law.

\subsection{Particle relabelling symmetry}
As we discussed in Section \ref{inverse map sec}, fluid equations in
general, and EPDiff in particular, are invariant under relabelling
of particles. In the context of the inverse map variables,
relabelling is accomplished by the action of the diffeomorphism
group $\Diff(\Omega)$ defined by
\[
\MM{l}\mapsto
\eta\circ\MM{l}\equiv\eta(\MM{l}), \qquad \eta\in\Diff(\Omega).
\]
The corresponding infinitesimal action
of the vector fields $\mathfrak{X}(\Omega)$ is then
\[
\MM{l} \mapsto \MM{\xi}\circ\MM{l}\equiv\MM{\xi}(\MM{l}), \qquad
\MM{\xi}\in\mathfrak{X}(\Omega),
\]
and the cotangent lift of this action is
\[
(\MM{\pi},\MM{l}) \mapsto
\left(-(\nabla\MM{\xi}(\MM{l}))^T\cdot\MM{\pi}
,\MM{\xi}(\MM{l})\right).
\]
To obtain the symmetry generator (\ref{X}), we extend the above
action to first derivatives as follows:
\begin{eqnarray*}
X &=& \xi_k(\MM{l})\pp{}{l_k} + (\xi_k(\MM{l}))_{,t}\pp{}{l_{k,t}}
+ (\xi_k(\MM{l}))_{,i}\pp{}{l_{k,i}} \\
 & & - \pi_k\pp{\xi_k(\MM{l})}{l_j}\pp{}{\pi_j}
- \left(\pi_k\pp{\xi_k(\MM{l})}{l_j}\right)_{,t}\pp{}{\pi_{j,t}} -
\left(\pi_k\pp{\xi_k(\MM{l})}{l_j}\right)_{,i}\pp{}{\pi_{j,i}}.
\end{eqnarray*}
The relabelling symmetries are variational, because
\[
XL = \pi_k(\xi_k(\MM{l}))_{,t} + \pi_ku_i(\xi_k(\MM{l}))_{,i} -
\pi_k\pp{\xi_k(\MM{l})}{l_j}\left(l_{j,t}+u_il_{j,i} \right)=0.
\]
Noether's theorem then gives the conservation law
\[
\left(\pi_k\xi_k(\MM{l})\right)_{,t}+\left(\pi_ku_j\xi_k(\MM{l})\right)_{,j}=0.
\]
A conservation law exists for each element $\MM{\xi}$ of
$\mathfrak{X}(\Omega)$, so particle relabelling generates an
infinite space of conservation laws.

\subsection{Circulation theorem}

To see how the particle relabelling conservation laws relate to
conservation of circulation, note that if $\rho$ is any density that
satisfies
\[
\rho_{,t} + (\rho u_j)_{,j} = 0,
\]
then
\[
\left( \frac{\pi_k\xi_k(\MM{l})}{\rho}\right)_{,t} + u_j \left(
\frac{\pi_k\xi_k(\MM{l})}{\rho}\right)_{,j}=0.
\]
If we pick a loop $C(t)$ which is advected with the flow, then
\[
\dd{}{t}\oint_{C(t)} \frac{\pi_k\xi_k(\MM{l})}{\rho}\rd x = 0.
\]
For a vector field $\MM{\xi}$ which is tangent to the loop at time $0$,
and satisfies $|\MM{\xi}|=1$ on the loop, then
\[
\MM{\xi}\rd x = (\nabla\MM{l})\cdot\rd{\MM{x}}
\]
for all times $t$, and one finds
\begin{equation}
\label{circulation}
\dd{}{t}\oint_{C(t)}\frac{\MM{\pi}\cdot(\nabla\MM{l})}{\rho}
\cdot\rd\MM{x}=0,
\end{equation}
The momentum formula (\ref{momentum map}) gives
\[
\dd{}{t}\oint_{C(t)}\frac{(1-\lambda^2\nabla^2)\MM{u}}{\rho}\cdot\rd\MM{x}=
\dd{}{t}\oint_{C(t)}\frac{\MM{m}}{\rho}\cdot\rd\MM{x}=0,
\]
which is the circulation theorem for EPDiff.


\section{Inverse map multisymplectic formulation for Euler-Poincar\'e
  equation with advected quantities}

\label{advected}

To extend this method to more general equations with advected quantities
is very simple: take the Lagrangian obtained from equation
(\ref{inverse map sec}) and add variables to represent
higher-order derivatives. For the sake of brevity we shall compute
one example, the incompressible Euler equations, and briefly discuss
the implications for the circulation theorem.

\subsection{Multisymplectic form of incompressible Euler equations}
We start with the reduced Lagrangian
\[
\ell[\MM{u},p,\rho] = \int_{\Omega}\frac{1}{2}\rho u_iu_i +
p(1-\rho)\rd V(\MM{x}),
\]
where $p$ is the pressure and $\rho$ is the relative density, and
add dynamical constraints to form the Lagrangian:
\[
L = \frac{1}{2}\rho u_iu_i + p(1-\rho) +
\pi_k\left(l_{k,t}+u_il_{k,i}\right) +\phi\left(\rho_{,t}+(\rho
u_i)_{,i}\right).
\]
This Lagrangian is already affine in the first-order
derivatives, so the Euler-Lagrange equations are automatically
multisymplectic in these variables:
\[
\begin{pmatrix}
0 & 0 & \pi_k\partial_i & 0 & -\rho\partial_i & 0 \\
0 & 0 & 0 & 0 & -\partial_t -u_i\partial_i & 0 \\
-\pi_k\partial_i & 0 & 0 & -\partial_t-u_i\partial_i & 0 & 0 \\
0 & 0 & \partial_t+u_i\partial_i & 0 & 0 & 0 \\
\rho\partial_i & \partial_t+u_i\partial_i  & 0 & 0 & 0 & 0 \\
0 & 0 & 0 & 0 & 0 & 0 \\
\end{pmatrix}
\begin{pmatrix}
u_i \\
\rho \\
l_k \\
\pi_k \\
\phi \\
p \\
\end{pmatrix}
= \nabla H,
\]
where the quantity
\[
H = -\left(\frac{1}{2}\rho u_iu_i + p(1-\rho)\right)
\]
is negative of the Hamiltonian density.

\subsection{Circulation theorem for advected quantities}
The conservation law for particle-relabelling follows exactly as in
Section \ref{inverse map EPDiff}, and we obtain equation
(\ref{circulation}) as before. The difference is that now the
momentum formula (momentum map) is
\[
\MM{m} = \pp{\ell}{\MM{u}} = -\pi_k\nabla l_k -\phi\diamond a
\]
and so one obtains
\[
\dd{}{t}\oint_{C(t)}\frac{\MM{m}}{\rho}\cdot\rd\MM{x}=
\oint_{C(t)}\frac{1}{\rho}\pp{\ell}{a}\diamond a\cdot \rd\MM{x}.
\]
For the incompressible Euler equations, $a$ is the relative density
$\rho$, so
\[
\pp{\ell}{a}\diamond a = \rho\nabla\pp{\ell}{\rho},
\]
which leads to the circulation theorem
\[
\dd{}{t}\oint_{C(t)}\frac{\MM{m}}{\rho}\cdot\rd\MM{x}=\oint_{C(t)}\nabla
\pp{\ell}{\rho}\cdot\rd\MM{x}=0.
\]
\section{A note on multisymplectic integrators}
\label{numerics} In this section we discuss briefly how to produce
multisymplectic numerical integrators, using the inverse map
formulation given in this paper. We note in particular that the
multisymplectic method will satisfy a discrete form of the
particle-relabelling symmetry and hence we will obtain a method that
has discrete conservation laws for $-\MM{\pi}\cdot\nabla\MM{l}$.

\subsection{Variational integrators}
A multisymplectic integrator for a PDE is a numerical method which
preserves a discrete conservation law for the two-form $\kappa$
given in equation (\ref{kappa}) (Bridges \& Reich, 2001). As
described in (Hydon, 2005), a discrete variational principle with a
Lagrangian that is affine in first-order differences automatically
leads to a set of difference equations which are multisymplectic.
This now makes it very simple to construct multisymplectic
integrators for fluid dynamics using the inverse map formulation:
one simply replaces the spatial and time integrals in the action
with numerical quadratures, replaces the first-order derivatives by
differences, and takes variations following the standard variational
integrator approach (Lew \emph{et al.}, 2003). Whilst the method
will preserve the discrete conservation law for the two-form
$\kappa$, the one-form quasi-conservation law will not be preserved
in general, and hence the other conservation laws will not be
exactly preserved.

\subsection{Discrete relabelling symmetry}
As $\MM{\pi}$ and $\MM{l}$ are still continuous in the discretised
equations, the multisymplectic integrator will have a discrete
particle-relabelling symmetry analogous to the one given in Section
\ref{inverse map EPDiff}, with the only difference being the
discretisation of the cotangent lift. Following the variational
integrator programme described in Lew \emph{et al.} (2003), the
discrete form of Noether's theorem will give rise to discrete
conservation laws for the multisymplectic method.

\subsection{Remapping labels}
If this approach is to be applied to numerical solutions with
intense vorticity then one needs to address the problem that
eventually the numerical discretisation of the labels $\MM{l}$ will
become very poor due to tangling, and hence the approximation to the
momentum
\begin{equation}
\label{mom} \MM{m} = -(\nabla\MM{l})^T\MM{\pi} = - \pi_k\nabla l_k
\,,
\end{equation}
will degrade with time. One possible approach would be to apply
discrete particle-relabelling, mapping the labels back to the
Eulerian grid in such a way that the momentum (\ref{mom}) stays
fixed. This transformation is exactly the relabelling given in
Section \ref{inverse map EPDiff}. Numerically, one could construct a
transformation (using a generating function for example) which
satisfies
\[
\MM{l} \mapsto \MM{X} + \mathcal{O}(\Delta x^p,\Delta t^p),
\qquad \MM{\pi}(\nabla\MM{l})\mapsto \MM{\pi}(\nabla\MM{l})
+ \mathcal{O}(\Delta x^p,\Delta t^p),
\]
where $p$ is the order of the method. For instance, one might use a
variational discretisation of the relabelling transformation, which
is generated by a symplectic vector field whose Hamiltonian is
$\MM{\pi}\cdot\MM{\xi}(\MM{l})=\pi_k\xi_k(\MM{l})$. In this way, one
may still retain some of the conservative properties of the method.

\section{Summary and Outlook}
\label{summary}
\subsection{Summary}

This paper describes a multisymplectic formulation of
Euler-Poincar\'e equations (which are, in essence, fluid dynamical
equations with a particle-relabelling symmetry). We have used the
inverse map to obtain a canonical variational principle, following
Holm and Kupershmidt (1983). As noted in Hydon (2005), a
multisymplectic formulation can be obtained by choosing variables
such that the Lagrangian at most linear in the first-order
derivatives, and contains no higher-order derivatives.  We have
shown how to construct the multisymplectic formulation for the
Euler-Poincar\'e equations for diffeomorphisms, using the example of
the EPDiff($H^1$) equations, and how to extend the method to the
Euler-Poincar\'e equations with advected quantities. These equations
encompass many fluid systems, including incompressible Euler,
shallow-water, Euler-alpha, Green-Naghdi, perfect complex fluids,
inviscid magnetohydrodynamics, \emph{etc.}

The techniques of Hydon (2005) have led to conservation laws for
these systems, including the usual multisymplectic conservation laws
for energy and momentum plus an infinite set of conservation laws
which arise from the particle-relabelling symmetry of fluid
dynamics.  We have highlighted the connection between these latter
conservation laws and Kelvin's circulation theorem, and showed that
multisymplectic integrators based on this formulation will have
discrete conservation laws associated with this symmetry.

\subsection{Outlook}

In the last section of this paper we have discussed the possibility
of developing multisymplectic integrators for fluids using this
framework. It is undoubtedly simple to construct such integrators,
but the issue of accuracy with time arises whenever the flow is
strongly mixing and numerical errors make the label field $\MM{l}$
very noisy. A discretisation of the relabelling map discussed in
this paper could provide a way to prevent this problem whilst
retaining some of the geometric properties of the method. These
ideas may aid the future development of integrators that have
conservation laws for vorticity and circulation, which are desirable
for numerical weather prediction and other applications.

In a different direction, we believe that multisymplectic
integrators would be especially apt for applications of EPDiff to
template-matching in computational anatomy (Holm \emph{et al}.,
2004). The matching problem is an initial-final value problem. In
such problems, space and time may be treated on an equal footing,
just as in the multisymplectic formulation.

\subsection{Acknowledgements}The work of DDH was partially supported
by the Royal Society of London Wolfson Award and the US Department of
Energy Office of Science ASCR.

\end{document}